\documentclass[a4paper,11pt]{article}
\usepackage{amsmath, amssymb}
\usepackage{doc, exscale, fontenc, latexsym, syntonly}
\usepackage{xypic}
\author{ A. Mutlu, T. Porter }

\title{\large{\bf Freeness Conditions for Crossed  Squares and Squared Complexes.}}

\newtheorem{prop}{Proposition}[section]
\newtheorem{thm}[prop]{Theorem}

\newtheorem{cor}[prop]{Corollary}

\newenvironment{pf}{{\bf Proof:}}{\hfill$\Box$\mbox{}}

\begin{document}

\maketitle
\begin{abstract}
Following Ellis, \cite{ellis2}, we investigate the notion of totally free
crossed square and related squared complexes.  It is shown how to interpret
the information in a free simplicial group given with a choice of CW-basis,
interms of the data for a totally free crossed square.  Results of Ellis then
apply to give a description in terms of tensor products of crossed modules. The 
paper ends with a purely algebraic derivation of a result of Brown and Loday. 

{ \noindent\bf A. M. S. Classification:}  18D35, 18G30, 18G50, 18G55, 55Q20, 55Q05.
\end{abstract}
\section*{Introduction}
Crossed squares were introduced by Loday and Guin-Walery in \cite{wl}. They
arose in various problems of relative algebraic K-theory.  Loday later showed
in \cite{loday} that these quite simple algebraic gadgets modelled all
homotopy 3-types. More generally his notion of cat$^n$-group and the related
crossed $n$-cubes of Ellis and Steiner were shown by  Loday to model all
connected $(n+1)$-types. The possibilities of calculation with these models was
enhanced by the development with R.Brown of a van Kampen type theorem for
these structures \cite{bl1}.

A link between simplicial groups and crossed $n$-cubes was used by Porter,
\cite{porter} to give an algebraic form of Loday's result and in particular to 
give a functor from the category of simplicial groups to that of crossed
$n$-cubes realising the equivalence.  

In 1993, Ellis \cite{ellis2} introduced a notion of free crossed square and
showed how to assign a free crossed square to a CW-complex. As there was an
established notion of free simplicial group, it seemed important to
investigate the extent to which the two notions of freeness are related. That
was the initial motivation for this paper. The two notions were intimately
related and moreover combining this with Ellis' alternative description of free crossed
squares in terms of the Brown-Loday non-abelian tensor product of groups and
coproducts of crossed modules, gives a new purely algebraic derivation of Brown 
and Loday's result describing the homotopy 3-type of the suspension of an
Eilenberg-Mac Lane space.  This success raises our hopes that this method of
attack can yield new results in higher dimensions.

\section{Preliminaries }
In this paper we will concentrate on the reduced case
and hence on simplicial groups rather than simplicial groupoids.  This is for
ease of exposition only and all the results do go through for simplicially
enriched groupoids. 

\textbf{Notation:} If $X$ is a set, $F(X)$ will denote the free group on $X$.
If $Y$ is a subset of $F(X)$, $\langle Y \rangle $ will denote the normal
subgroup generated by $Y$ within $F(X)$.

\subsection{Simplicial groups and groupoids}
Denoting the usual category of finite ordinals by $\Delta,$ we obtain for each 
$k\geq 0$, a subcategory $\Delta_{\leq k}$ determined by the objects $[j]$ of $\Delta$
with $j\leq k.$ A simplicial group is a functor from the opposite category
$\Delta^{op}$ to $\mathfrak{Grp};$ a $k$-truncated simplicial group is a functor 
from $\Delta^{op}_{\leq k}$ to $\mathfrak{Grp}.$ We will denote the  category of simplicial
groups by $\mathfrak{SimpGrp}$ and  the category of k-truncated 
simplicial groups by ${\mathfrak{Tr_kSimpGrp}}$. By a {\em k-truncation of a simplicial group}, we mean 
a $k$-truncated simplicial group $\mathfrak{tr_k}{\bf G}$ obtained by forgetting dimensions of
order $>k$ in a simplicial group {\bf G}, that is restricting {\bf G} to  $\Delta^{op}_{\leq k}$. This gives a truncation functor
$
\mathfrak{tr_k}:{\mathfrak{SimpGrp}}\longrightarrow {\mathfrak{Tr_kSimpGrp}} 
$
which admits a right adjoint
$
\mathfrak{ cosk_k}:\mathfrak{Tr_kSimpGrp}\longrightarrow \mathfrak{SimpGrp} 
$
called the {\em k-coskeleton functor}, and a left adjoint
$
\mathfrak{sk_k}:\mathfrak{Tr_kSimpGrp}\longrightarrow \mathfrak{SimpGrp,} 
$
called the {\em k-skeleton functor}. For explicit  constructions of these see
\cite{duskin}. We will say that a simplicial group $G$ is \emph{k-skeletal} if 
the natural morphism $\mathfrak{sk_k}G\rightarrow G$ is an isomorphism.

Recall that given a simplicial group {\bf G}, {\em the Moore complex} $(%
{ NG},\partial )$ {\em of} {\bf G} is the normal chain complex defined by 
$$
({ NG})_n=\bigcap_{i=0}^{n-1}\mbox{\rm Ker}d_i^n 
$$
with $\partial _n:NG_n\rightarrow NG_{n-1}$ induced from $d_n^n$ by
restriction. There is an alternative form of Moore complex given by the
convention of taking $$\bigcap^n_{i=1} \mbox{\rm Ker}d_i^n $$ and using $d_0$
instead of $d_n$ as the boundary.  One convention is used by Curtis 
\cite{curtis} (the $d_0$ convention) and the other by May \cite{may} (the $d_n$
convention). They lead to equivalent theories.

The {\em n$^{th}$ homotopy group} $\pi _n$({\bf G}) of {\bf G} is the $n
^{th}$ homology of the Moore complex of {\bf G}, i.e. 
$$
\begin{array}{rcl}
\pi _n({\bf G}) & \cong & H_n( 
{NG},\partial ) \\  & = & \bigcap\limits_{i=0}^n\mbox{\rm Ker}%
d_i^n/d_{n+1}^{n+1}(\bigcap\limits_{i=0}^n\mbox{\rm Ker}d_i^{n+1}). 
\end{array}
$$
We say that the Moore complex {\bf NG} of a simplicial group is of {\em length}
$k$ if $NG_n=1$ for all $n\geq k+1$, so that a Moore complex of length $k$ is also
of length $l$ for $l\geq k.$ For example, if ${\bf G}$ has Moore complex of length
1, then $(NG_1,NG_0,\partial_1)$ is a crossed module and conversely.  If $NG$
is of length 2, the corresponding Moore complex gives a 2-crossed module in
the sense of Conduch\'e, \cite{conduche}, cf. the companion
 paper to this, \cite{mp3}

\subsection{Free Simplicial Groups}
Recall from \cite{curtis} and \cite{kan2} the definitions of free simplicial
group  and of a $CW-basis$ for a free simplicial group.

{\bf Definition}

A simplicial group {\bf F} is called \emph{free} if\\
(a)\qquad $F_n$ is a free group with a given basis, for every integer $n\geq 0,$\\
(b)\qquad The bases are stable under all degeneracy operators, i.e., for every pair 
of integers $(i,n)$ with $0\leq i\leq n$ and every basic generator $x\in F_n$ the 
element $s_i(x)$ is a basic generator of $F_{n+1}.$

{\bf Definition}

Let ${\bf F}$ be a free simplicial group (as above). A subset $\mathfrak{F}\subset {\bf F}$ 
will be called a $CW-basis$ of ${\bf F}$ if \\
(a)\qquad $\mathfrak{F_n} = \mathfrak{F}\cap F_n$ freely generates 
$F_n$ for all $n\geq 0,$\\
(b)\qquad $\mathfrak{F}$ is closed under degeneracies, i.e. $x\in \mathfrak{F_n}$
implies $s_i(x)\in \mathfrak{F_{n+1}}$ for all $0\leq i\leq n,$\\
(c)\qquad if $x\in\mathfrak{F_n}$ is non-degenerate, then $d_i(x) = e_{n-1},$ the 
identity element of $F_n$, for all $0\leq i< n$.

As explained earlier, we have restricted attention so far to simplicial groups
and hence to connected homotopy types.  This is traditional but a bit
unnatural as all the results and definitions so far extend with little or no
trouble to simplicial groupoids in the sense of Dwyer and Kan \cite{D&K} and
hence to non-connected homotopy types.  It should be noted that such
simplicial groupoids have a fixed and constant simplicial set of objects and
so are not merely simplicial objects in the category of groupoids.  In this
context if $\mathbf{G}$ is a simplicial groupoid with set of objects $O$, the
natural form of the Moore complex $\mathbf{NG}$ is given by the same formula
as in the reduced case, interpreting Ker$ d^n_i$ as being the subgroupoid of
elements in $G_n$ whose $i^{th}$ face is an identity of $G_{n-1}$.
Of course if $n \geq 1$, the resulting $NG_n$ is a disjoint union of groups,
so $\mathbf{NG}$  is a disjoint union of the Moore complexes of the vertex
simplicial groups of $\mathbf{G}$ together with the groupoid $G_0$ providing
elements that allow conjugation between (some of) these vertex complexes
(cf. Ehlers and Porter \cite{ep}).

Crossed modules of, or over, groupoids are well known from the work of Brown
and Higgins.  The only changes from the definition for groups (cf.
\cite{loday}) is that one has to handle the conjugation operation slightly
more carefully:

A \emph{crossed module} is a morphism of groupoids  $\partial :
M\longrightarrow N$ where $N$ is a groupoid with object set $O$ say and
$M$ is a family of groups, $M = \{M(a) : a \in O\}$, together with  an action of
$N$ on $M$ satisfying (i) if $m \in M(a)$ and $n \in N(a,b) $ for $a,b,\in O$,
the result of $n$ acting on $m$ is ${}^nm \in M(b)$;
(ii) $\partial({}^nm)=n\partial(m)n^{-1}$ and (iii)
${}^{\partial(m)}{m'}=m{m'}m^{-1}$ for all $m,{m'}\in M$, $\ n\in N.$ For the
weaker notion in which condition (iii) is not required, the models are
called \emph{precrossed modules}.

The definition of a CW-basis likewise generalises with each  $\mathfrak{F}$  a
subgraph of the corresponding free simplicial groupoid.

\section{Crossed Squares and Simplicial Groups}
Although we will be mainly concerned with crossed squares in this paper, many
of the arguments either clearly apply or would seem to apply  in the more
general case of crossed $n$-cubes and $n$-cube complexes.  We therefore give
some background in this more general setting.  

Again although we give the definitions and results for groups, the
adaptation to handle groupoids over a fixed base is routine.

The following definition is due to Ellis and Steiner \cite{es}. Let  $<n>$
denote the set $\{1,...,n\}.$

{\bf Definition}

 A \emph{crossed $n$-cube} of groups is a family $ \{ \mathfrak{M}_A : A\subseteq <n>\}$  
of groups, together with homomorphisms $\mu _i:\mathfrak{M}_A\longrightarrow 
\mathfrak{M}_{A\setminus\{i\}}$ for $i\in <n>$ 
and functions 
$$h:\mathfrak{M}_A\times \mathfrak{M}_B\longrightarrow \mathfrak{M}_{A\cup B}$$for $A,B\subseteq <n>,$ such that if ${}^{a}b$ denotes $h(a,b)b$ for $a \in 
\mathfrak{M}_{A}$ and $b \in \mathfrak{M}_{B}$
with $A\subseteq B,$ then for all $a,{a'} \in \mathfrak{M}_{A}$ and $b,{b'} \in 
\mathfrak{M}_{B}, c \in \mathfrak{M}_{C}$ and 
$i,j \in <n>,$ the following  hold: 
$$
\begin{array}{ll}
1) & \mu _ia = a\ \quad 
\text{{\rm if}}\ i\not \in A, \\ 
2) & \mu _i\mu _ja = \mu _j\mu _ia, \\ 
3) & \mu _ih(a, ~b) = h(\mu _ia, ~\mu _ib), \\ 
4) & h(a, ~b) = h(\mu _ia, ~b)=h(a, ~\mu _ib)  \hfill 
\text{{\rm if}}\ i\in A\cap B, \\ 
5) & h(a, ~a^{\prime }) =\lbrack a,{~}a^{\prime }
\rbrack, \\ 
6) & h(a, ~b) = h{(b, ~a)}^{-1}, \\ 
7) & h(a, ~b)=1 \hfill\text{if $a = 1$ ~~\text{or}~~ $b = 1,$}\\
8) & h(aa^{\prime }, ~b)={}^{a} h(a^{\prime }, ~b)h(a, ~b), \\ 
9) & h(a,~bb^{\prime })= h(a,~b){~}{}^{b}h(a, ~b^{\prime }), \\ 
10) & {}^{a}h(b, ~c)=
h({}^{a}b,{~}{}^{a}c) \hfill \text{{\rm if }} A \subseteq B\cap C,\\ 
11)&{}^{a}h(h(a^{-1}, ~b), ~c)~{}^{c}h(h(c^{-1}, ~a), ~b)~{}^{b}h(h(b^{-1}, ~c), ~a) = 1. \\  
\end{array}
$$

{\em A morphism of crossed n-cubes} is defined in the obvious way: It is a
family of group homomorphisms, for $A\subseteq <n>,$ 
$
f_A:\mathfrak{M}_A\longrightarrow \mathfrak{M}_{A}^\prime
$
commuting with the $\mu _i$'s and $h$'s. We thus obtain a category of
crossed $n$-cubes which will be denoted by $\mathfrak{Crs^n},$ cf. Ellis and Steiner
\cite{es}. Again there is an obvious variant of this definition for groupoids
over a fixed set of objects, $O$.

{\bf Remark:}  Crossed squares, that is the case $n = 2$, were introduced by
Loday and Guin-Walery, \cite{wl}, but with an apparently different definition.
The two notions are however equivalent.

{\noindent{\bf Example 1:}}
For $n=1,$ a crossed 1-cube is the same as a crossed module.

 For $n=2,$ one
has a crossed 2-cube is a crossed square: 
$$
\diagram
\mathfrak{M}_{<2>} \dto_{\mu_1} \rto^{\mu_2} & \mathfrak{M}_{\{1\}} \dto^{\mu_1} \\
\mathfrak{M}_{\{2\}} \rto_{\mu_2} & \mathfrak{M}_{\emptyset}.
\enddiagram
$$
Each $\mu _i$ is a crossed module, as is $\mu _1\mu _2$. The $h$-functions give
actions and a function 
$$
h:\mathfrak{M}_{\{1\}}\times \mathfrak{M}_{\{2\}}\longrightarrow \mathfrak{M}_{<2>}. 
$$
The maps $\mu _2$  also define a map of crossed modules from
$(\mathfrak{M}_{<2>},\mathfrak{M}_{\{2\}} ,\mu_1)$ to $(\mathfrak{M}_{<1>},\mathfrak{M}_{\emptyset},\mu_1)$. In
fact a crossed square can be thought of as a crossed module in the category
of crossed modules.

{\noindent\bf Example 2:}
Let ${N_1}, {N_2}$ be normal subgroups  of a group $G$. 
The commutative square
diagram of inclusions;
$$
\diagram
{N_1} \cap {N_2}\dto_{~~~~~~~~~inc.} \rto^{inc.} & {N_2} 
\dto^{inc.} \\
{N_1} \rto_{inc.} & G   
\enddiagram 
$$
naturally comes together with actions of ${G}$ on ${N_1},{N_2}$ and 
${N_1}\cap {N_2}$ given by
conjugation and functions 
$$
\begin{array}{cccc}
h: & {N_A}\times {N_B}& \longrightarrow  & {N_A}\cap {N_B} = N_{A\cup B} \\  
& (n_1,n_2) & \longmapsto  & \lbrack n_1,~n_2\rbrack.
\end{array}
$$
That this is a crossed square is easily checked.

The following proposition is noted by the second author in  \cite{porter}.
\begin{prop}\label{alti} \cite{porter}
Let ${\bf G}$ be a simplicial group with simplicial normal subgroups 
${\bf N_1}$ and ${\bf N_2}.$ Then the square%
$$
\diagram
{\bf N_1}\cap {\bf N_2}\dto \rto &{\bf N_2} \dto \\
{\bf N_1} \rto &{\bf G}
\enddiagram
$$
induces a crossed square
$$
\diagram
\pi_0({\bf N_1}\cap {\bf N_2}) \dto \rto &\pi_0({\bf N_2})
\dto \\
\pi_0({\bf N_1})\rto &\pi_0({\bf G}).
\enddiagram
$$
\end{prop}
\begin{pf}
The $h$-function
$$
h:\pi _0({\bf {N}_1})\times \pi _0({\bf {N}_2})\longrightarrow 
\pi _0({\bf {N}_1}\cap {\bf N_2}) 
$$
is given by
$$
h({\overline{n_1}},~{\overline{n_2}})= {\overline{[n_1,~n_2]}} 
$$
for all $ {\overline{n_1}} \in \pi _0({\bf {N}_1}),\,{\overline{n_2}} \in
\pi _0({\bf {N}_2}).$ It is then simple, cf. \cite{porter}, to see that the
second diagram above is a 
crossed square. {~}\end{pf}\\
In fact up to isomorphism all crossed squares arise in this way,
cf. \cite{loday} and \cite{porter}. 

{\bf Example 3:}\label{z}
Let {\bf G} be a simplicial group. Let $\mathfrak{{M}}({\bf G},2)$ denote the 
following diagram
$$
\diagram
NG_2/\partial_3NG_3 \dto_{~~~~~\partial_2 '} \rto^{\qquad \partial_2} & NG_1\dto^{\mu} \\
\overline{NG_1} \rto_{~~~\mu'} & G_1 
\enddiagram
$$
Then this is the underlying square of a crossed square. The extra structure is given as follows: $NG_1=${\rm Ker}$d_0^1$ and $\overline{NG}_1=${\rm 
Ker}$d_1^1$.
Since $G_1$ acts on $NG_2/\partial _3NG_3,\ \overline{NG}_1$ and $NG_1,$
there are actions of $\overline{NG}_1$ on $NG_2/\partial _3NG_3$ and $NG_1$
via ${\mu'},$ and $NG_1$ acts on $NG_2/\partial _3NG_3$ and $\overline{NG}
_1$ via $\mu.$ Both $\mu $ and ${\mu'}$ are
inclusions, and all actions are given by conjugation. The $h$-map is 
$$
\begin{array}{ccl}
NG_1\times \overline{NG}_1 & \longrightarrow & NG_2/\partial _3NG_3 \\ 
(x,\overline{y}) & \longmapsto  & h(x,~y)= \lbrack s_1x,~s_1ys_0{y}^{-1} \rbrack  
\partial_3NG_3.
\end{array}
$$
Here $x$ and $y$ are in $NG_1$ as there is  a
bijection between $NG_1$ and $\overline{NG}_1.$  We leave the verification of
the axioms of a crossed square to the reader.\\
This example is clearly functorial and we denote by
$$
\diagram
\mathfrak{M}( - ,2)~:~\mathfrak{SimpGrp} \rto & \mathfrak{Crs^2},
\enddiagram
$$
the resulting functor.  This is the case $n = 2$ of a general construction 
of a crossed $n$-cube from a simplicial group given by the second author in
\cite{porter} based on some ideas of Loday.

{\noindent{\bf Examples 2 and 3 revisited:}}
Let $G$ be a group with normal subgroups ${N}_1, \ldots ,{N}_n$ of $G$. Let 
$$
\begin{array}{ccc}
\mathfrak{M}_A=\bigcap \{{N}_i:i\in A\} & \text{and} & \mathfrak{M}_\emptyset =${G}$
\end{array}
$$
with $A\subseteq <n>.$ For $i\in <n>,$   $\mathfrak{M}_A$ is a normal  
subgroup of $\mathfrak{M}_{A-\{i\}}$. Define 
$$
\mu _i:\mathfrak{M}_A\longrightarrow \mathfrak{M}_{A-\{i\}} 
$$
to be the inclusion. If $A,B\subseteq <n>$, then $\mathfrak{M}_{A\cup B}=\mathfrak{M}_A\cap \mathfrak{M}_B,$
let 
$$
\begin{array}{cccl}
h: & \mathfrak{M}_A\times\mathfrak{M}_B & \longrightarrow  & \mathfrak{M}_{A\cup B} \\  
& (a,b) & \longmapsto  & [{~}a,{~}b{~}]
\end{array}
$$
as $[\mathfrak{M}_A, \mathfrak{M}_B] \subseteq \mathfrak{M}_A\cap \mathfrak{M}_B,$ 
where $a\in \mathfrak{M}_A,\ b\in \mathfrak{M}_B.$ Then 
$$
\{\mathfrak{M}_A:\ A\subseteq <n>,\ \mu _i,\ h\} 
$$
is a crossed $n$-cube, called the {\em inclusion crossed n-cube} given by the
normal $n$-ad of groups $(G;\ {N}_1, \ldots, {N}_n).$

\begin{prop}\label{y}
Let $({\bf G};\ {N}_1, \ldots, {N}_n)$ be a simplicial normal $n$-ad of 
subgroups of groups and define for $A\subseteq <n>$%
$$
\mathfrak{M}_A=\pi _0(\bigcap\limits_{i\in A}{N}_i) 
$$
with homomorphisms $\mu _i:\mathfrak{M}_A\longrightarrow \mathfrak{M}_{A-\{i\}}$ and h-maps induced by
the corresponding maps in the simplicial inclusion crossed $n$-cube,
constructed by applying the previous example to each level. Then $\{\mathfrak{M}_A:\
A\subseteq <n>,\ \mu _i,\ h\}$ is a crossed $n$-cube.
\end{prop}
\begin{pf} See \cite{porter}.
\end{pf}\\

This  describes a functor, \cite{porter}, from the category of
simplicial groups to that of crossed $n$-cubes of groups.

\begin{thm}\label{t1}
If {\bf G} is a simplicial group, then the crossed $n$-cube $\mathfrak{M}${\rm (}%
{\bf G},$n${\rm )} is determined by:

(i) for $A\subseteq <n>,$%
$$
\mathfrak{M}({\bf G},n)_A=\frac{\bigcap_{j\in A}\text{{\rm Ker}} d_{j-1}^n}{
d_{n+1}^{n+1}(\text{{\rm Ker}} d_0^{n+1}\cap \{\bigcap_{j\in A}\text{{\rm Ker}
} d_j^{n+1}\})}; 
$$
(ii) the inclusion 
$$
\bigcap_{j\in A}\text{{\rm Ker}} d_{j-1}^n\longrightarrow \bigcap_{j\in
A-\{i\}}\text{{\rm Ker}} d_{j-1}^n 
$$
induces the morphism 
$$
\mu _i:\mathfrak{M}({\bf G},n)_A\longrightarrow \mathfrak{M}({\bf G},n)_{A-\{i\}}; 
$$

(iii) the functions, for $A,B\subseteq <n>,$ 
$$
h:\mathfrak{M}({\bf G},n)_A\times \mathfrak{M}({\bf G},n)_B\longrightarrow 
\mathfrak{M}({\bf G},n)_{A\cup B} 
$$
are given by 
$$
h(\bar x,\bar y) = \overline{[x,~y]}, 
$$
where an element of $\mathfrak{M}({\bf G},n)_A$ is denoted by $\bar{x}$ with $x\in
\bigcap_{j\in A}${\rm Ker}$d_{j-1}^n.$
\end{thm}\hfill$\Box$

Some simplification is possible, again see \cite{porter} for the details.

\begin{prop}\label{w}
If {\bf G} is a simplicial group, then

i) for $A\subseteq <n>,\ A\neq <n>,$ 
$$
\mathfrak{M}({\bf G},n)_A\cong \bigcap_{i\in A}\text{{\rm Ker}}d_{i-1}^{n-1} 
$$
so that in particular, $\mathfrak{M}({\bf G},n)_\emptyset \cong G_{n-1}$; in every
case the isomorphism is induced by $d_0,$

ii) if $A\neq <n>$ and $i\in <n>,$ 
$$
\mu _i:\mathfrak{M}({\bf G},n)_A\longrightarrow \mathfrak{M}({\bf G},n)_{A\setminus \{i\}} 
$$
is the inclusion of a normal simplicial subgroup,

iii) for $j\in <n>,$ 
$$
\mu _j:\mathfrak{M}({\bf G},n)_{<n>}\longrightarrow \bigcap_{i\neq j}\text{{\rm Ker}}%
d_i^{n+1} 
$$
is induced by $d_n.$
\end{prop}\hfill$\Box$

Expanding this data out for low values of $n$ gives:\\
\noindent 1) For $n=0$, 
$$
\begin{array}{rcl}
\mathfrak{M}({\bf G},0) & = & G_0/d_1(
\text{{\rm Ker}}d_0,) \\  & \cong  & \pi _0(
{\bf G}), \\  & = & H_0(N{\bf G).}
\end{array}
$$
2) For $n=1,$ $\mathfrak{M}({\bf G},1)$ is the crossed module 
$$
\mu_1 :\text{{\rm Ker}}d_0^1/d_2^2(NG_2)\longrightarrow G_1/d_2^2(\text{{\rm Ker}}d_0^2). 
$$

3) For $n=2,$ $\mathfrak{M}({\bf G},2)$ is 
$$
\diagram
\mbox{\rm Ker}d^{2}_{0} \cap \mbox{\rm Ker}d^{2}_{1} / 
d^{3}_{3} (\mbox{ \rm Ker}d^{3}_{0} 
\cap \mbox{\rm Ker}d^{3}_{1} \cap \mbox{\rm Ker}d^{3}_{2})\dto_{\quad\mu_1} 
\rto^{\qquad\quad \mu_2} & \mbox{\rm Ker}d^{2}_{0}  /  d^{3}_{3}(\mbox{\rm Ker}d^{3}_{0} 
\cap \mbox{\rm Ker}d^{3}_{1}) \dto^{\mu_1} \\
\mbox{\rm Ker}d^{2}_{1}  /  d^{3}_{3}(\mbox{\rm  Ker }d^{3}_{0} 
\cap \mbox{\rm Ker}d^{3}_{2}) \rto_{~~~~~\mu_2}  & G_2 / d^{3}_{3}
(\mbox{\rm Ker}d^{3}_{0}).
\enddiagram
$$
By Proposition \ref{w}, this is isomorphic to 
$$
\diagram
NG_2 / d^{3}_{3} (NG_3 ) \dto_{\mu_1} \rto^{\quad\mu_2}   & \mbox{\rm Ker }d^{1}_{0} \dto^{\mu_1} \\
\mbox{\rm Ker }d^{1}_{1} \rto_{\mu_2} &  G_1,   
\enddiagram
$$
that is
$$
\mathfrak{M}(\mathbf{G}, 2) \cong
\left ( \diagram
NG_2 / \partial_3 (NG_3 ) \dto \rto & \mbox{ Ker }d_{0} 
\dto \\
\mbox{ Ker }d_{1} \rto  &   G_1    
\enddiagram \right)
$$
is a crossed square. Here the $h$-map is 
$$h: \mbox{Ker}d_0^1\times\mbox{Ker}d_1^1\longrightarrow NG_2/d_3^3(NG_3)$$
given by $h(x,y) = [s_1x, ~s_1ys_0y^{-1}]~\partial_3NG_3$, as before.

Note if we consider the above crossed square as a vertical morphism
of crossed modules, we can take its kernel and cokernel within the category of crossed
modules. In the above, the morphisms in the top left hand corner are 
induced from $d_2$ so
$$
\mbox{Ker}\left ( 
\mu_1 : \frac{NG_2}{\partial_3NG_3}\longrightarrow  \mbox{Ker}d_1
\right) =  \frac{NG_2\cap \mbox{Ker}d_2}{\partial_3NG_3} 
 \cong  \pi_2({\bf G})
$$
whilst the other map labelled $\mu_1$ is an inclusion so has trivial kernel.
Hence the kernel of this morphism of crossed modules is
$$
\pi_2({\bf G})\longrightarrow 1.
$$
The image of $\mu_2$ is closed and normal in both the  groups 
on the bottom line and as $\mbox{Ker}d_0 =NG_1$ with the corresponding $\mbox{Im}\mu_1$
being $d_2NG_2,$ the cokernel is $NG_1/\partial_2NG_2,$ whilst 
$G_1/\mbox{Ker}d_0\cong G_0,$ i.e., the cokernel of $\mu_1$ is 
$\mathfrak{M}({\bf G}, 1).$

In fact of course $\mu_1$ is not only a morphism of crossed modules, it is 
a crossed module. This means that $\pi_2({\bf G})\longrightarrow 1$ is in 
some sense a $\mathfrak{M}({\bf G}, 1)$-module and that 
$\mathfrak{M}({\bf G}, 2)$ can be thought of as a crossed extension of  
$\mathfrak{M}({\bf G}, 1)$ by $\pi_2({\bf G}).$

\section{Free Crossed Squares}
\subsection{Definitions}
G. Ellis, \cite{ellis2}, in 1993 presented the notion of a free crossed
square. In this section, we recall his  definition  and give a construction of  free
crossed squares by using the second dimensional Peiffer elements and the $2$-skeleton of a
`step-by-step' construction of a free simplicial group with given $CW$-basis.
We firstly recall  the definition of a free crossed square on a pair of
functions $(f_2,f_3)$, as given by Ellis.  We will call these crossed squares
\emph{totally free}.

Let ${\bf B_1}, \ {\bf B_2}$ and ${\bf B_3}$ be sets. Take $F({\bf B_1})$ to
be the free group on ${\bf B_1}.$ 
Suppose given a function $f_2:{\bf B_2\longrightarrow }F({\bf B_1}).$ 
Let $\partial : M\longrightarrow F({\bf B_1})$ be the free pre-crossed
module on $f_2.$ Using the action of $F({\bf B_1})$ on $M$ we can form the
semi-direct product $M\rtimes F({\bf B_1}).$ The canonical inclusion $\mu: M
\longrightarrow M\rtimes F({\bf B_1})$ given by $m\mapsto(m,1)$
allows us to consider $M$ as a normal subgroup of $M\rtimes F({\bf B_1}).$ 
(Recall that any normal inclusion is a crossed module with action given by conjugation.)
There is a second normal subgroup of $M\rtimes F({\bf B_1})$ arising from $M,$ namely
$$
N = \{(m,\partial m^{-1}): m\in M\}\subset M\rtimes F({\bf B_1})
$$
with inclusion denoted ${\mu'} : N\longrightarrow M\rtimes F({\bf B_1}).$
For $m\in M$, we let ${m'}$ denote the element $(m^{-1},~\partial m)$ in $N.$

Assume given a function $f_3: {\bf B_3}\longrightarrow M,$ whose image lies in the kernel of the
homomorphism $\partial:M\longrightarrow F({\bf B_1}).$ There is then a corresponding function
${f_3'} : {\bf B_3}\longrightarrow N$ given by $y \mapsto (f_3(y), 1).$ 

{\bf Definition} \cite{ellis2}

A crossed square, 
$$\xymatrix{L\ar[r]^{\partial_2}\ar[d]_{\partial_2^\prime}& M\ar[d]^\mu \\
N\ar[r]_{\mu^\prime\hspace{5mm}}& M\rtimes F({\bf B_1}),}$$
is totally free  on the pair of functions $(f_2, f_3)$ if \\
(i) $(M,F({\bf B_1}), \partial)$ is the free pre-crossed module on $f_2$;\\
(ii)$\quad {\bf B_3}$ is a subset of $L$ with $f_3$ and ${f_3'}$ the 
restrictions of $\partial_2$ and ${\partial_2'}$ respectively; \\
(iii) for any crossed square$$\xymatrix{L^{\prime}\ar[r]^{\tau}\ar[d]_{\tau^\prime}& M\ar[d]^\mu \\
N\ar[r]_{\mu^\prime\hspace{5mm}}& M\rtimes F({\bf B_1}),}$$
and any function $\nu :{\bf B_3}\longrightarrow {L'}$ satisfying ${\tau}\nu = f_3,$
there is a unique morphism  $\Phi= (\phi, 1, 1, 1)$ of crossed squares: 
$$
\diagram
L \ddto_{\partial_2'} \drto^{\phi}  \rrto^{\partial_2} && {M} \xto'[1,0]_{\mu}[2,0] \drto^{=} \\
&{{L'}} \ddto\save\go[0,0];[1,0]:(0.5,-0.15)\drop{\tau^\prime}\restore \rrto^{\qquad\tau} && {M} \ddto^{\mu} \\
N \drto_{=} \xline'[0,1]^{\mu'}[0,2]|>\tip && M \rtimes {F({\bf B_1})} \drto^{=} \\
&N \rrto^{\mu'}  && M \rtimes F({\bf B_1})\\
\enddiagram
$$
such that $\phi{\nu'} =\nu,$ where ${\nu'}:{\bf B_3} \longrightarrow L$ is the 
inclusion.

We denote such a  totally free crossed square by $(L,M,N,M\rtimes
{F({\bf B_1})})$ omitting the structural morphisms from the notation when
there is no danger of confusion.

We know the free pre-crossed module on $f_2: {\bf B_2} \longrightarrow F({\bf B_1})$ is  $\partial : \langle{\bf B_2}\rangle \longrightarrow F({\bf B_1})$, where $\langle{\bf B_2}\rangle$ denotes the normal closure of ${\bf B_2}$ in the free 
group $F({\bf B_2}\cup s_0({\bf B_1}))$, 
so the function $f_3 : {\bf B_3}\longrightarrow M ~~( =\langle{\bf B_2}\rangle )$ 
is precisely the data $({\bf B_3}, f_3)$ for 2-dimensional construction data in the simplicial context, cf. \cite{mp2}. We thus need to recall the $2$-dimensional construction for a free simplicial group. This $2$-dimensional form
can be summarised by the diagram 
$${\bf \mathbb{F}}^{(2)}:
\diagram
...{~}F(s_1s_0({\bf B_1})\cup s_0({\bf B_2})
\cup s_1({\bf B_2})\cup {\bf B_3})
\rto<0.25ex> \rto<1ex> \rto<1.75ex>^{\qquad\qquad\quad d_0 ,d_1 ,d_2 } & 
F(s_0({\bf B_1})\cup {\bf B_2}) \lto<0.75ex> \lto<1.50ex>^{\qquad\qquad\quad s_1,s_0}
\rto<0.25ex> \rto<1ex>^{\qquad d_1, d_0} & F({\bf B_1}) \lto<0.75ex>^{\qquad s_0} 
\enddiagram $$
 with the simplicial morphisms given as in \cite{mp2}. 

\subsection{Free crossed squares exist.}
\begin{thm} 
A totally free crossed square $(L, M, N, M\rtimes {F({\bf B_1})})$ exists on the 2-dimensional construction data and is given by $\mathfrak{M}(\bf{F}^{(2)},2)$
where $\bf{F}^{(2)}$ is the 2-skeletal free simplicial group defined by the construction data.
\end{thm}
\begin{pf}
Suppose given the 2-dimensional construction data for a free simplicial
group, $\mathbf{F},$ which we will take as above as the data for a totally free crossed square. We will not assume detailed knowledge of \cite{mp2} so we start with $F({\bf B_1})$ and $f_2:{\bf B_2}\longrightarrow
F({\bf B_1})$ and form $M=\langle{\bf B_2}\rangle.$ This gives $\partial_1:
\langle{\bf B_2}\rangle \longrightarrow F(X_0)$ as the free pre-crossed module on $f_2.$ The semidirect product gives
$$
F(s_0({\bf B_1})\cup {\bf B_2})\cong M\rtimes F({\bf B_1})
$$
and we can identify this with $\bf{F}_1^{(2)}.$ This identification also
makes 
$$
M\cong \mbox{Ker}d_0^1
$$
for the $d_0^1$ of $\mathbf{F}^{(2)}.$
 
Next form $N =\{  (m,\partial m^{-1})\in M\rtimes F({\bf B_1}): m \in
M\}$. As $m\in\langle{\bf B_2}\rangle$, it is a product of conjugates of
elements of ${\bf B_2}$ and their inverses, so writing $m=\prod(m_{\alpha_i}) y_{\alpha_i}^{\varepsilon_i}(m_{\alpha_i})^{-1}$
for indices $\alpha_i$, and $\varepsilon_i =\pm
1$, we get $\partial m = \prod m_{\alpha_i}
t_{\alpha_i}^{\varepsilon_i}m_{\alpha}^{-1}$ 
where $t_i =f_2(y_i)$, which is also $ d_0^1(y_i)$.  Thus we can identify $N$ with $\langle \{ys_1d_0^1(y)^{-1}: y \in {\bf B_2}\}\rangle$, which is exactly Ker$d_1^1$.

Now 
$f_3:{\bf B_3}\rightarrow \mbox{\rm Ker}\partial_1 = \mbox{\rm Ker}
(\partial : NF^{(2)}_1\rightarrow NF^{(2)}_0)\subset  \langle {\bf B_2}\rangle.$
We know that this allows us to construct ${\bf F}_2^{(2)}$
and hence ${\bf F}_n^{(2)}$ for ~$n\geq 3$,
and in addition that taking 
$$L=NF^{(2)}_{2}/\partial _3(NF^{(2)}_3),$$
gives a crossed square
$$
\diagram
L\dto_{\partial'}\rto^{\partial} &  {M}\dto^{\mu}\\
N\rto_{\mu'} & F^{(2)}_1
\enddiagram
$$ 
which is  $\mathfrak{M}({\mathbf{F}}^{(2)}, 2).$ We claim this is the totally free crossed square on the construction data.

At this stage it is worth noting that there seems to be no    simple adjointness
statement between $\mathfrak{M}(-,2)$ and some functor that would  give  a quick proof of freeness. The problem is that $\mathfrak{M}(-,2)$
seems to be an adjoint only up to some sort of coherent homotopy. To 
avoid this difficulty we use a more combinatorial approach involving the higher dimension Peiffer elements and the explicit  description
of $L$.

In \cite{mp1}, we analysed in general the structure of groups of boundaries
such as $\partial _3 (NF_3^{(2)})$. There we showed that  $ NF_3^{(2)}$ is normally generated by elements of
the following forms:-\\
(i) For all $x\in NF^{(2)}_1,~y\in NF^{(2)}_2,$  
$$
\begin{array}{lcl}
f_{(1,0)(2)}(x , y) & = & [s_1s_0(x) , s_2(y)] [s_2(y) , s_2s_{0}(x)], \\ 
f_{(2,0)(1)}(x , y) & = & [s_2s_0(x) , s_1(y)] [s_1(y) , s_2s_1(x)] 
[s_2s_1(x) ,s_2(y)] [s_2(y) , s_2s_0(x)] ; 
\end{array}
$$
(ii) for all $ y \in NF^{(2)}_{2} , x \in NF^{(2)}_{1},$ 
$$
\begin{array}{lcl}
f_{(0)(2,1)}(x , y) & = & [s_0(x) , s_2s_1(y)] [s_2s_1(y) , s_1(x)] 
[s_2(x) , s_2s_1(y)], 
\end{array}
$$
and (iii) for all $x , y \in NF^{(2)}_2$, 
$$
\begin{array}{rcl}
f_{(0)(1)}(x , y) & = & [s_0(x) , s_1(y)] [s_1(y) , s_1(x)] [s_2(x) ,s_2(y)], \\ 
f_{(0)(2)}(x , y) & = & [s_0(x) , s_2(y)],  \\ 
f_{(1)(2)}(x , y) & = & [s_1(x) , s_2(y)] [s_2(y) , s_2(x)]. 
\end{array}
$$
Given our description of $NF^{(2)}$ in low dimensions, it is routine to
calculate normal generators  of the various groups involved here in terms of
 ${\bf B_1}$ and ${\bf B_2}$. We set $$Z=\{s_1(y)^{-1}s_0(y): y\in{\bf B_2}\}.$$
The above diagram can then be realised as
$$
\diagram
{J}\dto_{\partial_2'} \rto^{\partial_2} &{\langle{\bf B_2}\rangle} \dto^{\mu} & \\
\langle Z\rangle\rto_{{\mu'}\qquad} & \langle{\bf B_2}\rangle\rtimes F({\bf B_1}) 
\enddiagram .$$
Here $J$ is $(\langle s_1({\bf B_2})\cup {\bf B_3}\rangle \cap \langle
Z\cup{\bf B_3}\rangle )/P_2$,  $P_2$ being the second dimensional Peiffer normal subgroup, which is in fact 
just $\partial _3 (NF_3^{(2)})$, and which  is a subgroup of  
$\langle s_1({\bf B_2})\cup {\bf B_3}\rangle \cap\langle  Z\cup{\bf
B_3}\rangle$.

Given any crossed square $({L'},M, N, M\rtimes {F({\bf B_1})})$ and a function
$\nu: {\bf B_3}\longrightarrow {L'},$ there then exists a unique morphism
$$
\phi :(L,M,N,M\rtimes {F({\bf B_1})})\longrightarrow ({L'},M,N,
M\rtimes {F({\bf B_1})}) 
$$
given by%
$$
\phi ({y_i'}P_2)=\nu ({y_i'}) 
$$
such that $\phi{\nu'} =\nu.$ The existence of $\phi$ follows by using the
freeness property of the group $NF_2^{(2)}$ and then restricting to $\langle s_1({\bf B_2})\cup {\bf B_3}\rangle \cap\langle  Z\cup{\bf
B_3}\rangle$.  The normal generating 
elements of $P_2$ are then easily shown to have trivial image in $L'$ as that
group is part of the second crossed square.

Thus the diagram is the desired totally free crossed square on the
2-dimensional construction data. The crossed square properties of
$(L,M,N,M\rtimes {F({\bf B_1})}$ may be easily verified or derived from the  fact that this is exactly  $\mathfrak{M}({\bf{F}}^{(2)}, 2)$.
\end{pf} 

\medskip

{\bf Remark:}

At this stage, it is important to note that nowhere in the argument was use
made of the freeness of the 1-skeleton.  If $G$ is any 1-skeletal simplicial 
group and we form a new simplicial group
$H$ by adding in a set ${\bf B_3}$ of new generators in dimension 2, so that for 
instance, $H_2 = G_2 * F({\bf B_3})$, then we can use $M = NG_1 = {\rm Ker
  }d_0^{G,1}$ as before even though it need not be free.  The corresponding
$N$ is then isomorphic to ${\rm Ker }d_1^{G,1}$ with the bottom right hand
corner being $G_1$.  The `construction data' is now replaced by data for
killing some elements of $\pi_1(G)$, specified by $f_3 :{\bf B_3} \rightarrow
M$.  Although slightly at variance with the terminology used by Ellis,
\cite{ellis2}, we felt it  sensible to introduce
the term ``totally free crossed square'' for the type of free crossed square
constructed in the above theorem, using ``free crossed square'' for the more
general situation in which $(M,G,\partial)$ and $f_3$ are specified and no
requirement on 
 $(M,G,\partial)$ to be a free precrossed module is made.
\subsection{The $n$-type of the $k$-skeleton}
As in the other papers in this series, we will use the `step-by-step'
construction of a free simplicial group to observe the way in which the models 
react to the various steps of the construction.

In a `step-by-step' construction of a free simplicial group, there are
simplicial inclusions 
$$
{\bf{F}}^{(0)}\subseteq {\bf{F}}^{(1)}\subseteq 
{\bf{F}}^{(2)}\ldots 
$$
In general, considering the functor, $\mathfrak{M}(\quad ,n)$, from the category
of simplicial groups to that of crossed $n$-cubes, gives the corresponding morphisms 
$$
\mathfrak{M}({\bf{F}}^{(0)},\ n) \rightarrow \mathfrak{M}({\bf 
\bf{F}}^{(1)},\ n)
\rightarrow \mathfrak{M}({\bf{F}}^{(2)},\ n)\rightarrow  ... \rightarrow \mathfrak{M}({\bf{F}},\ n).
$$
We will investigate ${\bf \mathfrak{M}(\bf{F}}^{(i)}{\bf ,\ }n{\bf )}$,
for $n=0,1,2, $ and varying $i$.

Firstly look at $\mathfrak{M}({\bf{F}}^{(0)}{\bf ,\ }n{\bf ),\ }$ where the 0-skeleton ${\bf{F}}^{(0)}\,$ can be thought of as simplifying to 
$$\begin{array}{lccc}
{\bf{F}}^{(0)}: & \cdots \longrightarrow {F({\bf B_1})}\longrightarrow 
{F({\bf B_1})}\longrightarrow {F({\bf B_1})}
\end{array}
$$
with the $d_i^n=s_j^n=\ $identity homomorphism on ${F({\bf B_1})}$.

For $n=0,\,$ there is an equality

$$
\mathfrak{M}({\bf{F}}^{(0)}{\bf ,\ }0{\bf )=}F_0^{(0)}/d_1
(\text{Ker}d_0)= {F({\bf B_1})}, 
$$
and so $\mathfrak{M}({\bf{F}}^{(0)}{\bf ,\ }0)$ is just the free group
of 0-simplices of $\bf{F}$.

For $n=1$, $\mathfrak{M}({\bf{F}}^{(0)}{\bf ,}1{\bf )}$ is 
$
NF_1^{(0)}/\partial _2NF_2^{(0)}\longrightarrow F_0. 
$
 It is easy to show that $NF_1^{(0)}/\partial _2NF_2^{(0)}$ is trivial  and hence 
$$
\mathfrak{M}({\bf{F}}^{(0)}{\bf ,\ }1{\bf )\cong (}1\longrightarrow 
{F({\bf B_1})}). 
$$

For $n=2$, \ $\mathfrak{M}({\bf{F}}^{(0)}{\bf ,\ }2{\bf )}$ is the trivial
crossed square 
$$
\diagram
NF_2 / d^{3}_{3} (NF_3 )\ddto\rrto&&\mbox{\rm Ker }d^{1}_{0}\ddto   
&&1 \ddto\rrto&& 1\ddto \\
&&& = \\
\mbox{\rm Ker }d^{1}_{1}\rrto&&F_1&&1\rrto&&F({\bf B_1}).
\enddiagram
$$
Next look at $\mathfrak{M}({\bf{F}}^{(1)}{\bf ,\ }n{\bf )}$ and recall that the
1-skeleton{~} $\bf{F}^{(1)}$ is 
$$
\diagram
{\bf{F}}^{(1)} :\hspace{.5cm}
...{~}{F(s_1s_0({\bf B_1})\cup s_0({\bf B_2})\cup s_1({\bf B_2}))} 
\rto<0.25ex> \rto<1ex> 
\rto<1.75ex>^{\hspace{2.3cm} d_0 ,d_1 ,d_2 } & {F(s_0({\bf B_1})\cup {\bf B_2})} 
\lto<0.75ex> \lto<1.50ex>^{\hspace{2.4cm} s_1,s_0}
\rto<0.25ex>\rto<1ex>^{\qquad d_1, d_0} & {F({\bf B_1})} 
\lto<0.75ex>^{\qquad s_0}.
\enddiagram
$$

For $n=0$,  $\mathfrak{M}({\bf{F}}^{(1)},\ 0)$
is 
$
F_0^{(1)}/d_1(\text{Ker}d_0)\cong {F({\bf B_1})}/{\partial_1 NF_1}, 
$
which is $\pi _0({\bf{F}}^{(1)})\cong \pi _0(\bf{F}).$

For $n=1$, we have that 
$$
\begin{array}{rcl}
\mathfrak{M}({\bf\bf{F}}^{(1)}{\bf ,\ }1{\bf )} 
& = & (NF_1/\partial _2NF_2\longrightarrow
F_0), \\  
& = & 
\langle{\bf B_2}\rangle/P_1\longrightarrow {F({\bf B_1})},
\end{array}
$$
which is a free crossed module. In fact this is the free crossed module on the
(generalised) presentation $({\bf B_1};{\bf B_2}, f_2)$. As pointed out in
\cite{bh}, it is often convenient to generalise the notion of a presentation
$({\bf X}, {\bf R})$ with $${\bf R}\subset F(X)$$ to one  with the map ${\bf
  R} \rightarrow F(X)$ specified and not necessarily monic.  Thus if $f_2$ is
injective, this is just a presentation $\cal P$ of $\pi_1({\bf\bf{F}})$.
The kernel of this crossed module is then the module of identities of  $\cal
P$, again see \cite{bh}.

For $n=2$, $NF_2^{(1)} =\langle s_1({\bf B_2})\rangle \cap\langle Z\rangle$,
so $\mathfrak{M}({\bf{F}}^{(1)}{\bf ,\ }2{\bf )\ }$ 
simplifies to give (up to isomorphism), 
$$
\diagram
NF_2 / d^{3}_{3} (NF_3 )\ddto\rto&\mbox{\rm Ker }d^{1}_{0}\ddto 
&&J\ddto\rto&
\langle {\bf B_2}\rangle  \ddto \\
&& = \\
\mbox{\rm Ker }d^{1}_{1}\rto&G_1&& \langle Z\rangle \rto& 
F(s_0({\bf B_1})\cup {\bf B_2})
\enddiagram
$$
which is a crossed square with $J = 
(\langle s_1({\bf B_2})\rangle \cap\langle Z\rangle )/P_2.$

Next look at $\mathfrak{M}({\bf{F}}^{(2)}{\bf ,\ }n{\bf ).\ }$
Recall the 2-skeleton ${\bf{F}}^{(2)}$ is
$$
\diagram
{ \bf{F}}^{(2)}:\hspace{.5cm}
... {F(s_1s_0({\bf B_1})\cup s_0({\bf B_2})\cup s_1({\bf B_2})\cup{\bf B_3})} \rto<0.25ex> \rto<1ex> 
\rto<1.75ex>^{\qquad\hspace{1.8cm} d_0 ,d_1 ,d_2} & {F(s_0({\bf B_1})\cup {\bf B_2})} 
\lto<0.75ex> \lto<1.50ex>^{\qquad\hspace{1.8cm} s_1,s_0}
\rto<0.25ex> \rto<1ex>^{\qquad d_1, d_0} & {F({\bf B_1})} 
\lto<0.75ex>^{\qquad s_0}.
\enddiagram
$$
The following can be easily obtained by direct calculation :

 for $n=0,$ 
$$
\mathfrak{M}({ \bf{F}}^{(2)},0) = F_0/d_1(\text{Ker}d_0)\cong \pi _0({\bf{F}}^{(2)})= \mathfrak{M}({\bf{F}}^{(1)},0); 
$$

for $n=1,$ 
$$
\mathfrak{M}( \bf{F}^{(2)} ,1) \cong  \langle {\bf B_2}\rangle /P_1\longrightarrow {F({\bf B_1})}. 
$$
Finally, let $n=2.$ By an earlier result of this section, $\mathfrak{M}({\bf{F}}^{(2)}{\bf ,}2{\bf )\ }$ corresponds to the free crossed square,
$$
\diagram
NF_2 / d^{3}_{3} (NF_3 ) \ddto\rto &\mbox{\rm Ker }d^{1}_{0}\ddto 
&&J\ddto\rto &
\langle {\bf B_2}\rangle \ddto \\
&& = \\
\mbox{\rm Ker }d^{1}_{1}\rto&F_1&& \langle Z_2\rangle \rto 
& F(s_0({\bf B_1})\cup {\bf B_2})
\enddiagram
$$
where $J$ is now $(\langle s_1({\bf B_2})\cup{\bf B_3}\rangle \cap\langle
Z\cup{\bf B_3}\rangle )/P_2$ and $\langle Z_2\rangle $ is $\langle  Z\cup{\bf
  B_3}\rangle$, so this reduces to the earlier case if ${\bf B_3}$ is empty.
Thus we have the following relations 
$$
\mathfrak{M}({\bf{F}}^{(2)}{\bf ,\ }0{\bf )}=\mathfrak{M}({\bf{F}}^{(1)}{\bf ,\ }
0),
\qquad \mathfrak{M}({\bf{F}}^{(2)}{\bf ,\ }1{\bf )}=\mathfrak{M}({\bf{F}}^{(1)}
{\bf ,\ }1{\bf )} 
$$
but
$
\mathfrak{M}({\bf{F}}^{(2)}, 2)$ and $\mathfrak{M}({\bf{F}}^{(3)}, 2) $ need not be the same due to the additional influence of ${\bf B_3}$.
 Of course it is clear that, in general: 
$$
\begin{array}{ccccc}
 {\mathfrak{M}}({\bf{F}}^{(i)},\ n) & = & 
{\mathfrak{M}}({\bf{F}}^{(i+1)},\ n) & \text{if} & i\geq n+1. 
\end{array}
$$

\section{Squared Complexes}
The authors and Z. Arvasi have  defined $n$-crossed complexes  in \cite{zat}.
In this paper, we will only need the case $n= 2$, which had already been
defined by Ellis in \cite{ellis2}.   We shall follow him in calling these   
{\em squared complexes}. A squared complex consists of a diagram of group homomorphisms
$$
\diagram
& & & & N\drto^{\mu'}\\
\ldots \rto& C_4\rto^{\partial_4}&C_3\rto^{\partial_3}& 
L\urto^{\lambda'}\drto_{\lambda}&& P\\
& & & &M\urto_{\mu}  
\enddiagram
$$
together with actions of $P$ on $L, N, M$ and $C_i$ for $i\geq 3,$ and a function 
$h : M\times N\longrightarrow L.$ The following axioms need to be satisfied.\\
(i) The square $\left(\spreaddiagramrows{-1.2pc} \spreaddiagramcolumns{-1.2pc}
\def\objectstyle{\ssize} \def\labelstyle{\ssize}
\diagram
L \dto_{{\lambda'}} \rto^{\lambda}  & N \dto^{\mu} \\
{M} \rto_{\mu'}& P
\enddiagram\right)$  is a crossed square; \\ 
(ii) The group $C_n$ is abelian for $n \geq 3$; \\
(iii) The boundary homomorphisms satisfy $\partial_n\partial_{n+1} = 1$ for
$n \geq 3,$ and  $\partial_3(C_3)$ lies in the intersection 
$\mbox{ker}\lambda\cap \mbox{ker}{\lambda'};$\\
(iv)  The action of $P$ on $C_n$ for $n \geq 3$ is such that ${\mu}{M}$ and
${\mu'}N$ act trivially. Thus each $C_n$ is a $\pi_0$-module with $\pi_0 =
P/{\mu}{M}{\mu'}N$; \\
(v) The homomorphisms $\partial_n$ are $\pi_0$-module homomorphisms for $n \geq 3.$

This last condition does make sense  since the axioms for crossed squares imply that 
$\mbox{ker}{\mu'}\cap \mbox{ker}{\mu}$ is a $\pi_0$-module.

A morphism of squared complexes $$\Phi: (C_\ast,\left(\spreaddiagramrows{-1.2pc} 
\spreaddiagramcolumns{-1.2pc}
\def\objectstyle{\ssize} \def\labelstyle{\ssize}
\diagram
L \dto_{{\lambda'}} \rto^{\lambda}  & N \dto^{\mu} \\
{M} \rto_{\mu'}& P
\enddiagram\right) )\longrightarrow ({C_\ast'},\left(\spreaddiagramrows{-1.2pc} 
\spreaddiagramcolumns{-1.2pc}
\def\objectstyle{\ssize} \def\labelstyle{\ssize}
\diagram
{L'} \dto_{{\lambda'}} \rto^{\lambda}  & {N'} \dto^{\mu} \\
{M'} \rto_{\mu'}& {P'}
\enddiagram\right))$$
consists of a morphism of crossed squares $(\Phi_{L}, 
\Phi_{N}, \Phi_{M}, \Phi_{P})$, together with a family of equivariant homomorphisms $\Phi_n$ for $n \geq 3$ satisfying
$\Phi_{L}\partial_3 = {\partial'}_3\Phi_{L}$ and 
$\Phi_{n-1}\partial_n = {\partial'}_n\Phi_n$ for $n \geq 4.$ 
There is clearly a category $\mathfrak{SqComp}$ of squared complexes. This
exists in both group and groupoid based versions.

By a {\em (totally) free squared complex}, we will mean one in which the crossed square
is (totally) free,  and in which each $C_n$ is free as a $\pi_0$-module  for $i\geq 3.$

\begin{prop}

There is  a functor 
$${\mathcal C}(\quad ,2)
: \mathfrak{SimpGrp}\longrightarrow \mathfrak{SqComp}
$$
such that free simplicial groups are sent to totally free squared complexes.
\end{prop}

\begin{pf}

Let ${\bf G}$ be a simplicial group or groupoid. We will define a squared
complex ${\mathcal C}({\bf G},2)$ by specifying ${\mathcal C}({\bf G},2)_A$ for each
$A\subseteq <2>$ and for $n \geq 3$,  ${\mathcal C}({\bf G},2)_n$. As usual, (cf. the 
other papers in this series, \cite{mp, mp1, mp2,mp3}), we will
denote by $D_n$ the subgroup or subgroupoid of $NG_n$ generated by the
degenerate elements.

For $A\subset <2>,$ we define
$${\mathcal C}({\bf G},2)_A =
\mathfrak{M}(\mathfrak{sk_{2}}{\bf G,}2)_A = \frac{\cap\{\mbox{Ker}d_{i}^{2} : i \in A
\}}{d_3(\mbox{Ker}d^3_0 \cap \bigcap\{\mbox{Ker}d_{i+1}^{3} : i \in A \}
\cap D_3)}.
$$
We do not need to define $\mu_i$ and the $h$-maps relative to these groups as they
are already defined in the crossed square $\mathfrak{M}(\mathfrak{sk_{2}}{\bf
  G,}2)$.

For $n \geq 3$,  we set $${\mathcal C}({\bf G},2)_n = \frac{NG_n}{(NG_n\cap
  D_n)d_{n+1}(NG_{n+1}\cap D_{n+1})}.$$
As this is part of the crossed complex associated to ${\bf G}$, we can take the
structure maps to be those of that crossed complex, cf. \cite{ep,mp2}.  The terms are all modules 
over the corresponding $\pi_0$ as is easily checked.  The final missing piece, 
$\partial_3$, of the structure is induced by the differential $\partial_3$ of $NG$.

The axioms for a squared complex can now be verified using the known results
for crossed squares and for crossed complexes with a direct verification of
those axioms relating to the interaction of the two parts of the structure,
much as in \cite{ep} and \cite{mp2}.

Now suppose the simplicial group is free. The proof above of the freeness of
$\mathfrak{M}(\mathfrak{sk_{2}}{\bf G,}2)$ together with the freeness of the
crossed complex of a free simplicial group, \cite{mp2}, now completes the proof.
\end{pf}

\medskip

Suppose that $\rho$ is a general squared complex.
The {\em homotopy~groups} $\pi_{n}(\rho),$ $n\geq 0$ of $\rho$ are defined cf. \cite{ellis2}, to be the
homology groups of the complex
$$
\diagram
\ldots\rto^{\partial_5~}&C_4\rto^{\partial_4~~}&C_3\rto^{\partial_2~}&L\rto^{\partial_2~\quad}
&M\rtimes N\rto^{\quad\partial_1}&P\rto&1
\enddiagram
$$ 
with $\partial_2(l)=({\lambda'}l^{-1}, \lambda l)$ and $\partial_1(m,n) = \mu(m){\mu'}(n).$
The axioms of a crossed square guarantee that $\partial_2$ and $\partial_1$ are homomorphisms with
$\partial_3(C_3)$ normal in $\mbox{Ker}(\partial_2),~ \partial_2(L)$ normal
in $\mbox{Ker}(\partial_1),$ and $\partial_1(M\rtimes N)$ normal in $P$.

\begin{prop}
The homotopy groups of ${\mathcal C}({\bf G},2)$ are isomorphic to those of
${\bf G}$ itself.
\end{prop}
\begin{pf}

Again this is a consequence of well-known results on the two parts of the structure.
\end{pf}

\section{Alternative Descriptions of Freeness.}
In the context of CW-complexes, Ellis, \cite{ellis2} gave a neat description
of the top group $L$ in a (totally) free crossed square derived from that data.  A
simplicial group with a given CW-basis is the algebraic analogue of a
CW-complex so one would expect a similar result to hold in that setting.
Ellis uses the generalised van Kampen theorem of Brown and Loday, \cite{bl1}.
In the algebraic setting no such tool is available, but in fact its use is not 
needed.

Ellis' description  is in terms of tensor products and coproducts.  For
completeness we recall the background definitions of these constructions.

\subsection{Tensor Products}

Suppose that $\mu: M\to P$ and $\nu: N\to P$ are crossed modules over $P.$ The groups 
$M$ and $N$ act on each other, and themselves, via the action of $P.$ The tensor product 
$M\otimes N$ is the group generated by the symbols $m\otimes n$ for $m\in M$, $n\in N$ subject to 
the relations
$$mm'\otimes n = ({}^{m}{m'}\otimes {}^{m}n)(m\otimes n),$$
$$m\otimes n{n'} = (m\otimes n)({}^{n}{m}\otimes {}^{n}{n'}),$$
for $m,{m'}\in M, \ n,{n'}\in N.$ There are homomorphisms $\lambda: M\otimes N \to M, \
{\lambda'}: M\otimes N \to N$ defined on generators by $\lambda(m\otimes n)= m({}^{n}m)^{-1}$
and ${\lambda'}(m\otimes n)= ({}^{m}n)n^{-1}.$ The group $P$ acts on $M\otimes N$ by 
${}^{p}(m\otimes n) = ({}^{p}m\otimes {}^{p}n),$ and there is a function $h: M\times N
\to M\otimes N,$ $(m,n)\longmapsto m\otimes n.$ In \cite{bl1}, it is  verified  that this structure
gives a crossed square
$$
\diagram
M\otimes N \rto^{\lambda} \dto_{\lambda'}  & N\dto^{\nu} \\
M \rto_{\mu} & P  \\
\enddiagram
$$          
with the universal property of extending the corner $$
\diagram
 & N\dto^{\nu} \\
M \rto_{\mu} & P  \\
\enddiagram.
$$          
\subsection{Coproducts}
Let $(M, P,\partial_1), (N,P,\partial_2)$ be $P$-crossed modules. Then $N$ acts on $M,$ and $M$ 
acts on $N,$ via the given actions of $P.$ Let $M\rtimes N$ denote the semidirect product with 
the multiplication given by
$$
(m,n)({m'},{n'})=(m{m'},{~}  {}^{m'}n{n'})
$$
and injections
$$
\begin{array}{cc}
{i'}: M\to M\rtimes N \qquad \mbox{and} \qquad {j'} : N\to M\rtimes N \\
\quad m\longmapsto (m,1)\qquad{~}\qquad\qquad\qquad    n\longmapsto(1,n).
\end{array}
$$
We define the pre-crossed module
$$
\begin{array}{cc}
\underline{\delta}: M\rtimes N\to P \\
(m,n) \longmapsto \partial_1(m)\partial_2(n). 
\end{array}
$$
Let $\{M, N\}$ be the subgroup of $M\rtimes N$ generated by the elements of the form
$$
(m{}^nm^{-1},n{}^mn^{-1})
$$
for all $m\in M$, $n\in N$,  thus we are able to form the quotient group
$M\rtimes N/\{M,N\}$ and obtain an induced morphism
$$
\partial: M\rtimes N/\{M,N\}\to P
$$ 
given by 
$$
\partial(m,n)\{M,N\} = \partial_1(m)\partial_2(n).
$$
Let $q: M\rtimes N\to M\rtimes N/\{M,N\}$ be projection and let $i = q{i'}, \ j= q{j'}.$ Then
$M\circ N = (M\rtimes N)/\{M,N\}$ with the morphisms $i,j,$ is {\em a
  coproduct} of $(M,P,\partial_1)$ and $(N,P,\partial_2)$ 
in the category of $P$-crossed modules.
\begin{prop}~\cite{ellis2}
Let $(L, M, \bar{M}, M\rtimes F)$ be a (totally) free crossed square on the
$2$-dimensional  construction data or on functions $(f_2, f_3)$ as described above. Let $\partial: C\to M\rtimes F$ be the free crossed module on the function ${\bf B_3}\to M\rtimes F$
given by $y\longmapsto (f_3y, 1).$ From the crossed module $M\otimes \bar{M} \to
M\rtimes F$, then $L$ is isomorphic to the coproduct
$(M\otimes \bar{M})\circ C$
factored by  the relations
$$
\begin{array}{cc}
1)\quad i(\partial c\otimes\bar{m}) = j(c)j({}^{\bar{m}}c^{-1}) \\
2)\quad i(m\otimes\partial c) = j({}^{m}c)j(c^{-1})
\end{array}
$$
for $c \in C, \ m\in M ~\mbox{and}~ \bar{m}\in\bar{M}.$

The homomorphisms $L\to M$, $ L\to \bar{M}$ are given by the homomorphisms
$$
\lambda:M\otimes \bar{M}\to M\qquad\mbox{and}\qquad {\lambda'}: M\otimes \bar{M}\to \bar{M}
$$
and $\partial : C\to M\cap\bar{M}.$
The $h$-map of the crossed square is given by 
$$
h(m, \bar{n}) = i(m\otimes\bar{n})
$$
for $m, n \in M.$
\end{prop}
\begin{pf} This comes by direct verification using the universal properties
of tensors and coproducts. \end{pf}

\medskip

\textbf{Remark: }  For future applications it is again important to note that the result 
is not dependent on the crossed square being \emph{totally} free, although
this is the form proved and used by Ellis, \cite{ellis2}.  If $M\rightarrow
F$ is any pre-crossed module, one can form the `corner'
$$
\diagram
 & M\dto \\
\bar{M} \rto & M\rtimes F,  \\
\enddiagram$$         
complete it to a crossed square via $M\otimes \bar{M}$ and then add in ${\bf
  B_3}\to M$. Nowhere does this use freeness of $M\rightarrow
F$.
\begin{cor}
Let $\bf{G}^{(1)}$ be the $1$-skeleton of a simplicial group.
Then in  the free crossed square $\mathfrak{M}(\bf{G}^{(1)}, 2)$ described above, 
$$
NG_2^{(1)}/\partial_3NG_3^{(1)}\cong \mbox{Ker}d_1^{1}\otimes\mbox{Ker}d_0^{1}.
$$
\end{cor}
\begin{pf}
This is clear from the previous proposition.
\end{pf}

\medskip

\textbf{Remark: }

If we set $M = \mbox{Ker}d^1_0 = NG^{(1)}_1$, then the identification given by the Corollary gives
$$NG^{(1)}_2/\partial_3NG^{(1)}_3 \cong M\otimes \bar{M}.$$ This uses the fact
that $\mbox{Ker}d^1_0 $ and $\mbox{Ker}d^1_1$ are linked via the map sending
$m$ to $ms_0d_1m^{-1}$ for $m \in \mbox{Ker}d^1_0 $.  The $h$-map $h : M \times \bar{M} \rightarrow
NG^{(1)}_2/d^3_3NG^{(1)}_3$ is $h(x,y) = [s_1x,s_1ys_0y^{-
1}]d^3_3NG^{(1)}_3$, but this is also $h(x,y) = x\otimes y$. Thus 
$$x\otimes y = [s_1x,s_1ys_0y^{-1}]d^3_3NG^{(1)}_3$$
under the identification via the isomorphism of 5.2.

This explains the `mysterious' formula of \cite{mp} in the discussion
before Proposition 4.6 of that paper.

\subsection{Applications to 2-crossed complexes.}

Of course there are similar results for free squared complexes.  What
is less obvious is the way in which these results can be applied to the
situation that we studied in our earlier paper, \cite{mp3}.  There we
considered the alternative model for 3-types given by Conduch\'e's 2-crossed
modules and also looked at the corresponding 2-crossed complexes. We will
not repeat all  that discussion here but note the definition: 

{\bf Definition:}\\
A 2-crossed complex of group(oid)s is a sequence of group(oid)s 
$$C:\hspace{1cm} \ldots \rightarrow C_n
\stackrel{\partial_n}{\rightarrow}C_{n-1}\rightarrow\ldots C_2
\stackrel{\partial_2}{\rightarrow}C_1\stackrel{\partial_1}{\rightarrow}C_0$$
in which\\
(i) ~$C_n$ is abelian for $n\geq 3$;\\
(ii)~$C_0$ acts on $C_n$, $n\geq 1$, the action of $\partial C_1$ being
trivial on $C_n$ for $n\geq 3$;\\
(iii)~ each $\partial_n$ is a $C_0$-group(oid) homomorphism and
$\partial_i\partial_{i+1} =1$ for all $i\geq 1$;\\
and\\
(iv)~ $C_2
\stackrel{\partial_2}{\rightarrow}C_1\stackrel{\partial_1}{\rightarrow}C_0$
is a 2-crossed module.

\medskip

We refer the reader to \cite{conduche} or \cite{mp3} for the exact meaning of
2-crossed module.

Given a simplicial group or groupoid, ${\bf  G}$, define
$$C_n = \left\{\begin{array}{ll}
        NG_n & \mbox{ \rm for } n = 0,1\\
        NG_2/d_3(NG_3 \cap D_3)& \mbox{ \rm for } n = 2\\
        NG_n/(NG_n \cap D_n)d_{n+1}(NG_{n+1} \cap D_{n+1})& 
\mbox{ \rm for } n \geq 3
\end{array}\right.
$$
with $\partial_n$ induced by the differential of $\bf{NG}$.  Note that the
bottom three terms (for $n = $ 0, 1, and  2) form a 2-crossed module
considered in \cite{conduche} or \cite{mp3} and that for $ n \geq 3$, the groups are all $ \pi_0(G)$-modules, since in these dimensions $C_n$ is
the same as the corresponding crossed complex term (cf. Ehlers and
Porter \cite{ep} for instance).

\begin{prop}\cite{mp3}

With the above structure $(C_n, \partial_n)$ is a 2-crossed
complex, which will be denoted $C(\mathbf{G})$.\hfill$\Box$
\end{prop}

Here we note in particular that the  term $C_2$ is $ NG_2/d_3(NG_3 \cap D_3)$
and so is the same as ${\mathcal C}({\bf G},2)_{<2>}$.  Thus 
if ${\bf G}$ is a simplicial group, we obtain {\em gratis}:
\begin{cor}
Let $\bf{G}^{(1)}$ be the $1$-skeleton of a simplicial group.
The 2-crossed complex of $\bf{G}^{(1)}$ satisfies  $$C(\mathbf{G}^{(1)})_2\cong
 \mbox{Ker}d_1^{1}\otimes\mbox{Ker}d_0^{1}.
$$\hfill$\Box$
\end{cor}
We also get in general a description of $C(\bf{G}^{(2)})_2$ as a quotient
of the form $( \mbox{Ker}d_1^{1} \otimes\mbox{Ker}d_0^{1}\circ C)/\sim$ where as in
Proposition 5.1, this $C$ is a free crossed module on the `new cells' in
dimension 2.

\subsection{The suspension of a $K(\pi,1)$.}

As was mentioned in \cite{mp2}, Brown and Loday used their
generalised van Kampen Theorem, \cite{bl1}, to calculate
$\pi_3\Sigma K(\pi,1)$ for $\pi$ a group, as the kernel of the commutator
map from $\pi \otimes \pi$ to $\pi$.  Jie Wu, (\cite{wu} Theorem 5.9),
 for any group $\pi$ and set of generators $\{x_\alpha |
\alpha \in J\}$ for $\pi$, gives a presentation of 
$\pi_n\Sigma K(\pi,1)$ in terms of higher commutators, but does not
manage to get the Brown-Loday result explicitly although his result is
clearly linked to theirs.

Wu's methods use a study of simplicial groups and a construction he
ascribes to Carlsson, \cite{carlsson}.  This gives a simplicial group
$F^\pi(S^1)$ that has $\pi_{n+2}\Sigma K(\pi,1)\cong \Omega \Sigma K(\pi,1) \cong \pi^{n+1}  
F^\pi(S^1)$.  As we pointed out in \cite{mp}, $F^\pi(S^1)$ is a
pointed analogue of the `tensorisation' of $K(\pi,0)$, the constant
simplicial group on $\pi$, with the simplicial circle $S^1$.  In general
if $G$ is a simplicial group and $K$ a pointed simplicial set,
$G\bar{\wedge}K$ will denote the simplicial group with group of $n$-
simplices given by 
$$
\coprod\limits_{x\in K_n}(G_n)_{x}/(G_n)_{\ast}.
$$
If $x \in K_n$, we denote the $x$-indexed copy of $g \in G_n$ within $(G\bar{\wedge}K)_n$
by $g\bar{\wedge} x$. The face and degeneracy maps of $G\bar{\wedge}K$  are
induced by the componentwise application of the corresponding morphisms of $G$
and $K$
$$d_i(g\bar{\wedge}x) = d_i^Gg\bar{\wedge}d_i^Kx,$$
$$s_i(g\bar{\wedge}x) = s_i^Gg\bar{\wedge}s_i^Kx.$$
Of course if $d_i^Kx = \ast$ then $d_i(g\bar{\wedge}x) = 1$.

The case of interest to us is $G = K(\pi,0)$, $K = S^1$ and we will adopt the
notation for simplices in $S^1$ used by us in \cite{mp}.  We write $S^1_0 = \{ 
\ast\}$ and will take $\ast$ to denote the corresponding degenerate
$n$-simplex basing $S^1_n$ in all dimensions; $S^1_1 = \{\sigma, \ast\}$,
$S^1_2 = \{x_0, x_1, \ast\}$, where $x_0 = s_1\sigma$, $x_1 = s_0\sigma$ and
in general $S_{n+1}^1 = \{x_0, \ldots, x_n, \ast\}$, where $x_i = s_n\ldots
s_{i+1}s_{i-1} \ldots s_0\sigma$, $0\leq i \leq n$.

We write $G = K(\pi,0)$ for simplicity and will usually make no distinction
between simplices in different dimensions unless confusion might arise.  We
have

 $(G\bar{\wedge}S^1)_0 = 1$, \quad the trivial group,

$(G\bar{\wedge}S^1)_1 \cong \pi,$

$(G\bar{\wedge}S^1)_2 \cong \pi \ast \pi$, \qquad the free product of two copies 
of $\pi$, and so on.  \\
The group $(G\bar{\wedge}S^1)_n$ is a free product of
$n$-copies of $\pi$, $\coprod\{(\pi)_x : x \in S^1_n\setminus \{\ast\}\}$, and 
writing as above $g\bar{\wedge}x$ for the $x$-indexed copy of $g\in \pi$ in
this, we note that $(g\bar{\wedge}x)(g^\prime\bar{\wedge}x) =
(gg^\prime\bar{\wedge}x)$ for $g$, $g^\prime \in \pi$.
As $g\bar{\wedge}x_i^{(n+1)} = s_n(g\bar{\wedge}x_i^{(n)})$ holds in all
dimensions, $n\geq 2$ and for all $0 \leq i \leq n$, it is clear that
$N(G\bar{\wedge}S^1)_n = D_n$, that is, it is generated by degenerate elements 
in all dimensions $n \geq 2$,  we can therefore apply Corollary 5.2.  As $N(G\bar{\wedge}S^1)_0$ is trivial, $\mbox{Ker}d_0^1 = \mbox{Ker}d_1^1 =(G\bar{\wedge}S^1)_1 \cong \pi$, so we get:

For $H = G\bar{\wedge}S^1$,$$NH_2/\partial_3NH_3 \cong \pi \otimes \pi.$$
We have by \cite{porter} that the algebraic 2-type of $H$ is completely
modelled by the crossed square $\mathfrak{M}(H,2)$, that is by 
$$
\diagram
\pi\otimes\pi \rto^{\mu_2}\dto_{\mu_1}&\pi\dto^=\\
\pi\rto_=&\pi,
\enddiagram
$$
where $\mu_1$ and $\mu_2$ are the commutator maps.

As a consequence we have:
\begin{cor}
The 3-type of $\Sigma K(\pi,1)$ is completely specified by the above crossed
square.  In particular there is an isomorphism
$$\pi_3(\Sigma K(\pi,1)) \cong \mbox{Ker}(\mu : \pi \otimes \pi \rightarrow
\pi).$$
\hfill$\Box$
\end{cor}
This result was first found by Brown and Loday \cite{bl1}.  Their proof was an 
illustration of the use of their generalised van Kampen Theorem.  Jie Wu,
\cite{wu}, gives some methods that shed light on the higher homotopy groups,
but although they yield a description of $\pi_4$, they  do not analyse the
4-type itself.
The model $G\bar{\wedge}S^1$ is 1-skeletal and one might  expect that a triple
tensor $\pi \otimes\pi \otimes\pi$ may be involved in any model of its 4-type.
  Of course $\mathfrak{M}(H,3)$
gives a complete model, but the individual terms involved in that model are
not as easy to analyse as in $\mathfrak{M}(H,2)$. An amalgam of Wu's methods and
the methods developed in the earlier papers of this series,
\cite{mp,mp1,mp2,mp3}, might  provide insight into this. This problem
is not of itself that important, but it does seem to provide an excellent
testbed for the development of methods to aid in calculation with low
dimensional algebraic models of homotopy types.


{\noindent  A. Mutlu}\hspace{5.7cm} \ \ \ T. Porter  \\      
{\it Department of Mathematics, \ \ \hspace{2.5cm}School of Mathematics,\\
Faculty of Science, \hspace{3.6cm} \ \ \ \ \ University of Wales, Bangor, \\
University of Celal Bayar, \hspace{3cm}\ Gwynedd, LL57 1UT, UK.\\
Manisa, Turkey. \hspace{4.6cm} UK\\
e-Mail: amutlu@spil.bayar.edu.tr \hspace{2.1cm}e-Mail: t.porter@bangor.ac.uk

\end{document}